\documentclass[a4paper,fleqn]{cas-dc}

\usepackage{amssymb}
\usepackage[numbers]{natbib}
\usepackage{hhline}
\usepackage{threeparttable}
\usepackage{tabularx}
\usepackage{multirow}
\usepackage{xcolor}
\usepackage{booktabs}
\usepackage{multicol}
\usepackage[ruled,vlined]{algorithm2e}
\usepackage[figuresright]{rotating} 
\usepackage{subfigure}
\usepackage{graphicx}
\usepackage{colortbl}
\usepackage{color}
\usepackage[numbered]{bookmark}
\usepackage{hyperref}
\hypersetup{
        bookmarksnumbered=true,
        bookmarksopen=true,
        bookmarksopenlevel=5,
        pdfstartview=Fit,
        pdfpagemode=UseOutlines
    }
\usepackage{array}

\usepackage[switch]{lineno}

\def\tsc#1{\csdef{#1}{\textsc{\lowercase{#1}}\xspace}}
\tsc{WGM}
\tsc{QE}
\tsc{EP}
\tsc{PMS}
\tsc{BEC}
\tsc{DE}

\begin{document}
\let\WriteBookmarks\relax
\def\floatpagepagefraction{1}
\def\textpagefraction{.001}
\shorttitle{A Knowledge-driven MA for the EEDHFSSP}
\shortauthors{Yunbao Xu et~al.}

\title [mode = title]{A Knowledge-driven Memetic Algorithm for the Energy-efficient Distributed Homogeneous Flow Shop Scheduling Problem}                      



\author[1]{Yunbao Xu}
\ead{08036@hnie.edu.cn}
\address[1]{School of Management, Hunan Institute of Engineering, Xiangtan 411104, China}

\author[2]{Xuemei Jiang}
\ead{xuemeijane@126.com}
\address[2]{Key Laboratory of Collaborative Intelligence Systems, Ministry of Education, Xidian University, Xi’an 710071, China}

\author[1]{Jun Li}
\ead{jli@hnie.edu.cn}

\author[2]{Lining Xing}
\ead{lnxing@xidian.edu.cn}
\cormark[1]

\author[3]{Yanjie Song}
\ead{songyj_2017@163.com}
\address[3]{School of Reliability and Systems Engineering, Beihang University, 100191, China.}

%
%
%
%
%
%
%
%
%

\cortext[cor1]{Corresponding author}
%

\begin{abstract}
The reduction of carbon emissions in the manufacturing industry holds significant importance in achieving the national "double carbon" target. Ensuring energy efficiency is a crucial factor to be incorporated into future generation manufacturing systems. In this study, energy consumption is considered in the distributed homogeneous flow shop scheduling problem (DHFSSP). A knowledge-driven memetic algorithm (KDMA) is proposed to address the energy-efficient DHFSSP (EEDHFSSP). KDMA incorporates a collaborative initialization strategy to generate high-quality initial populations. Furthermore, several algorithmic improvements including update strategy, local search strategy, and carbon reduction strategy are employed to improve the search performance of the algorithm. The effectiveness of KDMA in solving EEDHFSSP is verified through extensive simulation experiments. It is evident that KDMA outperforms many state-of-the-art algorithms across various evaluation aspects.
\end{abstract}




\maketitle

\section{Introduction}

At present, China is in a critical phase of accelerating the construction of a new development paradigm and deepening comprehensive reforms. Building a green and low-carbon industrial system and development pattern is an intrinsic requirement for promoting high-quality economic development, practicing new development concepts, and fostering harmonious coexistence between humans and nature \cite{yan2024distributed}. Therefore, driving carbon emission reduction in the manufacturing industry is of key significance for achieving China’s dual carbon goals. Production scheduling, as a vital component of the manufacturing system, directly affects the efficiency and competitiveness of enterprises \cite{LR2023TASE}. Researching and applying efficient optimization techniques and scheduling methods is crucial for achieving energy conservation and emission reduction, reducing production costs, and enhancing the optimality of production systems \cite{yao2019multi}. It is also central to improving production efficiency, economic benefits, and considering environmental impacts \cite{lu2020energy}. With the rapid development of economic globalization, an increasing number of manufacturing enterprises are expanding their production models into distributed environments and establishing multiple factories in different geographical locations to enhance production efficiency and cope with intensifying market competition \cite{pan2024bi}. Carbon emissions during the production scheduling process vary depending on the different processing jobs and technologies \cite{luo2024q}. It is evident that different production scheduling plans significantly affect the total carbon emissions of the entire production process. Therefore, investigating the distributed flow shop scheduling problem (DFSSP) based on carbon emissions not only has profound academic significance but also holds considerable practical application value.

In the DFSSP, numerous challenges arise from the realistic constraints or unique stages of different flow shop types. Existing research has considered constraints related to flow time, including fuzzy processing times, random processing times, setup times, and transportation times. Constraints related to production flow shops include no-wait, no-idle, blocking, limited buffer, and batch flow. Other constraints considered in distributed flow shop scheduling include job reentrant \cite{RIFAI201642}, unrelated machines \cite{doi:10.1080/00207543.2019.1598596}, and heterogeneous production flow shops \cite{MENG2021100804, LI2020106638}. The methods for solving distributed flow shop scheduling problems are mainly divided into exact and approximate methods. Exact methods are only suitable for small-scale problems. When the problem size becomes too large, they may result in excessively long computation times or become unsolvable. Approximate methods include heuristic algorithms, metaheuristic algorithms, and hybrid algorithms \cite{liu2024tri,xie2023hybrid,wang2020bi}. Heuristic algorithms can produce scheduling plans in a relatively short time, but the quality of their solutions is difficult to guarantee. In contrast, metaheuristic and hybrid algorithms can produce better approximate solutions within an acceptable time frame \cite{HZLG202206001}. Regarding scheduling objectives, those related to completion time and machine workload are the most evaluated. Energy consumption and low-carbon-related objectives are increasingly receiving attention and can be considered as one of multiple objectives, optimized simultaneously with traditional objectives. In DFSSP research, the use of real-world cases is relatively uncommon. Most instances are expanded from benchmarks of classic flow shop and job shop scheduling problems. For adaptive searching, various strategies are employed to improve their local and global search performance \cite{9409755}. 

To Addressing the issue of carbon emissions in flow shop scheduling, Ai et al. proposed a novel neighborhood search, based on the problem's characteristics, that combines memory and global exchange. This method is designed to solve the hybrid flow shop scheduling problem with the goal of reducing carbon emissions \cite{ai2017novel}. Zhong et al. considered both economic and environmental factors, setting the minimization of the make-span and total carbon emissions as optimization goals, and introduced a hybrid cuckoo search algorithm for solving the multi-objective permutation flow shop scheduling problem \cite{ZGJX201822004}. Geng et al. proposed an improved multi-objective hybrid memetic algorithm (MA) designed to solve the green reentrant hybrid flow shop scheduling problem, with objectives to minimize the make-span, total energy consumption costs, and carbon emissions \cite{ZGJX202012013}. Jubiz-Diaz et al. developed two multi-objective models integrating packaging size and production scheduling in flexible flow shop system, minimizing both costs and carbon emissions while determining the packaging dimensions and production plans for each finished product, and introduced a Pareto-based hybrid genetic algorithm \cite{doi:10.1080/00207543.2019.1566650}. Wu et al. studied the multi-objective flexible flow shop batch scheduling problem considering variable processing and handling times, with the goals of minimizing make-span and carbon emissions, and proposed a hybrid non-dominated sorting genetic algorithm with variable local search \cite{doi:10.1504/IJAAC.2020.110071}. Gu et al. also took into account economic and environmental factors, researching the multi-objective permutation flow shop scheduling problem with objectives of minimizing the make-span and total carbon emissions, and proposed a hybrid cuckoo search algorithm \cite{doi:10.1177/16878140211023603}. Saber et al. aimed to minimize total delay and carbon emissions, presenting a multi-objective decomposition-based heuristic algorithm based on job insertion and a multi-objective variable neighborhood search (VNS) algorithm for solving the permutation flow shop scheduling problem \cite{paperid:1087151}. Fernandez-Viagas et al., focusing on green permutation flow shop scheduling with variable processing times, considered both make-span and total energy consumption cost, proposing an iterative local search algorithm based on critical paths \cite{10.1016/j.cie.2022.108276}.

Addressing the issue of carbon emissions in DFSSP, Dong et al. established a two-stage reentrant hybrid flow shop bi-level scheduling model with the optimization objectives of minimizing the make-span, total carbon emissions, and total energy consumption costs. They proposed an improved hybrid salp swarm and NSGA-III algorithm for solving this problem \cite{zhang2022matrix}. Zhang et al. introduced a matrix cube-based distributed estimation algorithm aimed at solving the energy-efficient distributed assembly and permutation flow shop scheduling problem, which minimizes both make-span and total carbon emissions \cite{SJES5C54B49C5235AAA98149FBB04E2F1B27}. Schulz et al. studied distributed permutation flow shop scheduling with the goal of minimizing the make-span and carbon emissions. They used an adaptive bisection method for optimizing small instances and analyzing problem characteristics, and for practical situations, they developed a new multi-objective iterative greedy algorithm \cite{SJESF8B5C0D487BBA57A7CABC24E674AB4CC}. Shao et al., aiming to minimize total tardiness, total production cost, and total carbon emissions, established a distributed heterogeneous hybrid flow shop scheduling model considering energy and labor, and proposed a network MA for solution \cite{KZYC201601001}.

\textcolor[rgb]{0,0,0}{The performance of evolutionary algorithms in various types of flow shop scheduling problems is commendable \cite{zhao2023knowledge,zhang2024multiobjective}. Many researchers also design corresponding improvement strategies for specific problems \cite{song2024reinforcement}. In order to solve EEDHFSSP, we further improve the algorithm based on the MA framework according to the problem characteristics. The main contributions of this paper are as follows.}

\textcolor[rgb]{0,0,0}{1. A variant of DHFSSP that takes into account carbon emissions is investigated. accordingly, a multi-objective mathematical model is constructed to optimize both makespan and energy consumption. A number of factors such as processing sequence, equipment capacity, etc. are considered in the constraints.}

\textcolor[rgb]{0,0,0}{2. A knowledge-driven memetic algorithm (KDMA) is proposed. KDMA designs collaborative initialization strategy, updating strategy, local search strategy based on key factories. in addition, a carbon reduction strategy for energy consumption is used in the algorithm. carbon reduction strategy is used in the algorithm to save energy.}

\textcolor[rgb]{0,0,0}{3. Simulation experiments are conducted to test the effectiveness of KDMA for solving EEDHFSSP. KDMA outperforms other comparison algorithms in terms of diversity of solutions, convergence, and overall performance. In addition, the effectiveness of the strategies in the algorithm is also verified in the experiments.}

\textcolor[rgb]{0,0,0}{The remainder of this paper is organized as follows. Section \ref{Distributed Flow Shop Scheduling Problem} introduces the DFFSP. Section \ref{model} introduces the model of the EEDHFSSP. The KDMA is proposed in Section \ref{KDMA}. Section \ref{Experiment} verifies the effect of the proposed algorithm. Section \ref{Conclusion} summarizes the content and analyzes possible directions for further research in the future.}

\section{Distributed Flow Shop Scheduling Problem}
\label{Distributed Flow Shop Scheduling Problem}
\begin{figure}[htp]
\centering
\includegraphics[width=0.45\textwidth]{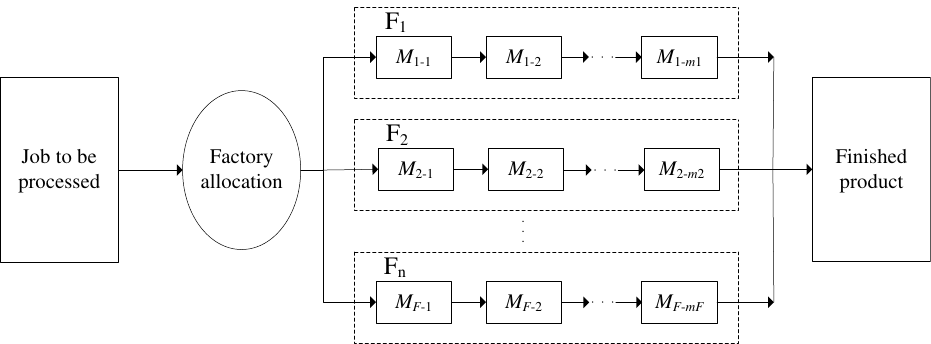}
\caption{Process flow diagram of distributed flow shop scheduling}
\label{Process flow diagram of distributed flow shop scheduling}
\end{figure}

The distributed shop scheduling problem (DSSP) refers to the collaborative production between companies or different factories, jointly processing and manufacturing products. The core of this problem lies in allocating jobs to different factories and arranging the processing sequence of these jobs within each factory, to optimize one or multiple scheduling indices \cite{lei2021improved}. The DSSP is an extension of the flow shop scheduling problem and is classified as an NP-hard problem.
DSSPs are mainly divided into distributed parallel machine scheduling \cite{sym14020204,PMID:35942738}, distributed flow shop scheduling problem (DFSSP) \cite{iceis17}, and distributed job shop scheduling problem \cite{csahman2021discrete}. Among these, the DFSSP is widely prevalent in the manufacturing industry. It involves allocating  \(n\) jobs to  \(F\) flow shops for processing, where each factory contains \(m\)  machines. The processing sequence of each job is fixed. In other words, every job goes through the machines in the same order. The process flow diagram of the DFSSP is shown in Figure.\ref{Process flow diagram of distributed flow shop scheduling}.

The DFSSP can be categorized based on the processing technology into distributed permutation flow shop scheduling and distributed flexible flow shop scheduling \cite{han2021distributed}. Distributed permutation flow shop is a specific case of distributed flow shops, requiring that the processing sequence of jobs on each machine is also identical. Distributed flexible flow shop refers to flow shops in factories where at least one operation involves multiple parallel processing machines. Furthermore, DFSSPs can be classified according to the type of factory into distributed homogeneous flow shop scheduling problem (DHFSSP) and distributed heterogeneous flow shop scheduling problem. In homogeneous flow shops, the dispersed multiple flow shops have identical processing capabilities, meaning the type and number of machines, the machines' age, processing precision, and types of auxiliary materials consumed are all the same. These can be seen as multiple factories established according to a standard model, where each job can be processed in any factory with the same processing conditions. In contrast, in heterogeneous flow shops, there are differences in the processing capabilities of each factory, such as the age of processing equipment, processing precision, and machine power. These differences lead to variations in process parameters for jobs in different factories, and thus different scheduling schemes significantly impact the scheduling indices.

\section{Mathematical Model}
\label{model}

\subsection{Symbols and Variables}

The parameters and variables used in this study are defined in Table \ref{Parameter definitions for the EEDHFSSP model}.

\begin{table}[htbp]
	\centering
	\caption{Parameter definitions for the EEDHFSSP model}
	\label{Parameter definitions for the EEDHFSSP model}
	\begin{tabularx}{0.45\textwidth}{m{1.5cm}m{5cm}}
		\toprule[1.5pt]
		Parameter & Definition \\
		\midrule[0.75pt]
		$n$ & Number of jobs \\
		$m$ & Number of machines in each factory \\
		$F$ & Number of factories \\
		$i$ &  Index for jobs, $i\in \left \{ 1,2,..., n \right \}$  \\
		$j$ & Index for machines, $j\in \left \{ 1,2,...,m \right \}$\\
		$f$ & Index for Factories, $f\in \left \{ 1,2,...,F \right \}$\\
		$k$ & Index of the job position in a given sequence, $k\in \left \{ 1,2,...,n \right \} $\\
		$P_{i,j} $ & Processing time of job $i$  on machine $j$  \\
		$C_{k,j}$ & Completion time of the job at position $k$ on machine $j$   \\
		$M$ & An infinitely large positive number \\
		$C_{max}$ & makespan \\
		$PP_{i,j} $ & Power consumption of job $i$  on machine $j$ during processing \\
		$SP$ & Power consumption of the machine when idle \\  
		$\varepsilon _{elec} $ & Carbon emission coefficient for electricity \\  
		$\varepsilon _{au} $ & Carbon emission coefficient for auxiliary materials\\  
		$CE$ & Total carbon emissions \\  
		$CE_{run}  $ & Carbon emissions generated during job processing \\  
		$CE_{idle}  $ & Carbon emissions generated during machine idle time \\  
		$CE_{au}  $ & Carbon emissions generated from the consumption of auxiliary materials for machines \\  
		\bottomrule[1.5pt]                            
	\end{tabularx}
\end{table}

\subsection{Problem Description}

The distributed homogeneous flow shop scheduling problem based on carbon emission (EEDHFSSP) is described as follows: $n$ jobs are allocated to $F$ identical flow shops for processing, each job comprising $M$ processing operations. All shops have identical flow-line processing equipment, with each flow line comprising $M$ processing machines. The processing time for each job is the same in different factories. In the flow shops, the sequence of processing on the $M$ machines is the same for each job, and the order of processing for jobs on each machine in the same flow shop is identical. To summarize, EEDHFSSP is the process of determining the start time of each operation of a job to complete the production of all the products as soon as possible, taking into account the energy consumption. 

\subsection{Model}

This problem is studied based on the following assumptions.

(1) All jobs are independent and can start processing at zero time.

(2) Jobs can be processed in any factory. Once a job is assigned to a factory, all its processing must be completed within that factory.

(3) Each machine can process only one job at a time, and each job can be processed by only one machine at any given moment.

(4) Once a job begins processing on a designated machine, it cannot be pre-empted or interrupted.

(5) The total carbon emissions during the entire production scheduling process include three parts: carbon emissions generated during job processing, emissions during machine idle times, and emissions from the consumption of auxiliary materials (such as lubricants and coolants for machines).

(6) Machine set-up times and transfer times between job operations are negligible, and situations involving machine failure or insufficient memory are not considered.

In our work, the optimization objectives of EEDHFSSP are to minimize the makespan and total carbon emissions.

\textbf{Objective function:}
\begin{equation}
	\min  \left \{ C_{max} ,CE \right \}  
\end{equation}

In the EEDHFSSP model, let  $\Pi$  be the final arrangement of  $n$ jobs in $F$  factories, with the arrangement of jobs in the $f$-th factory denoted as $\pi ^{f}$. Given that the number of jobs allocated to the $f$-th factory is $ n _{f}$ , we have $n_{1} +  n_{2} +...+n_{f} = n$ , and $\pi ^{f} = \left \{ \pi ^{f}\left ( 1 \right )  ,..., \pi ^{f}\left ( n  _{f} \right )  \right \} $. The final arrangement $Pi$  contains $F$ arrangements, denoted by $\Pi=\left \{ \pi ^{1},\pi ^{2},...,\pi ^{f} \right \} $ . The method for calculating $C_{max}\left ( \Pi \right )$ consists of the following five steps:

\textbf{Step 1:}\ $C_{\pi ^{f}\left ( 1 \right ),1}  =p_{\pi ^{f}\left ( 1 \right ),1}, f\in \left \{ 1,2,...,F \right \}$ 

\textbf{Step 2:}\ $C_{\pi ^{f}\left ( i \right ),1}  =C_{\pi ^{f}\left ( i-1 \right ),1}+p_{\pi ^{f}\left ( i \right ),1},f\in \left \{ 1,2,...,F \right \},j\in \left \{ 2,3,...,m \right \}$

\textbf{Step 3:}\ $ C_{\pi ^{f}\left ( 1 \right ),j}  =C_{\pi ^{f}\left ( 1 \right ),j-1}+p_{\pi ^{f}\left ( 1 \right ),j},f\in \left \{ 1,2,...,F \right \},j\in \left \{ 2,3,...,m \right \}$

\textbf{Step 4:}\ $ C_{\pi ^{f}\left ( i \right ),j}  =\max \left \{  C_{\pi ^{f}\left ( i-1 \right ),j},C_{\pi ^{f}\left ( i \right ),j-1}\right \}+p_{\pi ^{f}\left ( i \right ),j}  ,f\in \left \{ 1,2,...,F \right \},i\in \left \{ 2,3,...,n_{f}  \right \},j\in \left \{ 2,3,...,m \right \}$

\textbf{Step 5:}\ $ C_{max} \left ( \Pi \right ) =\max \left (C _{n_{f},m }  \right )$

Therefore, in EEDHFSSP, the makespan is defined by the maximum completion time across all sub-arrangements in the final arrangement $\Pi$ , specifically, it is determined by the maximum completion time of job $ n _{f}$  on the last machine in the $f$-th factory, denoted by $ C_{max} \left ( \Pi \right ) =\max \left (C _{max }\left ( \pi_{f} \right )   \right )  =\max \left ( C_{n_{f},m} \right ) $ .

In the production scheduling process, carbon emissions originate from four sources: emissions during job processing, emissions when machines are idle, emissions from the transportation of jobs between machines, and emissions from the consumption of auxiliary materials during processing. In DFSSPs, the scheduling sequence exerts a minimal influence on emissions associated with the transportation of jobs between machines. Therefore, we only consider carbon emissions produced during job processing ($CE_{run}$), emissions during machine idle time ($CE_{idle}$), and emissions from the consumption of auxiliary materials during processing ($CE_{au}$). The total carbon emissions (CE) for EEDHFSSP can be represented as $CE=CE_{run}+CE_{idle}+CE_{au}$, with the calculation formulae for the three types of emissions as follows:

(1) Carbon emissions during job processing ($CE_{run}$):
\begin{equation}
	\setlength	\abovedisplayskip {3pt} 
	CE_{run}=\sum_{i=1}^{n} \sum_{k=1}^{n}\left ( x_{i,k}\cdot p_{i,j}\cdot PP_{i,j} \right ) \cdot \varepsilon_{elec}
	\setlength	\belowdisplayskip {3pt} 
\end{equation}

The carbon emissions during job processing equal the product of the actual processing time of all jobs, the rated energy consumption during processing, and the carbon emission coefficient for electricity. Here,  $p_{i,j}$ denotes the processing time of job $i$  on machine $j$. $PP_{i,j}$ represents the unit energy consumption of job $i$  during processing on machine $j$.  $x_{i,j}$ is a decision variable. $\varepsilon_{elec}$  is the electricity emission coefficient.

(2) Carbon emissions during machine idle time ($CE_{idle}$):
\begin{equation}
	\setlength	\abovedisplayskip {3pt} 
	CE_{idle}=\sum_{f}^{E} \left [ C\left ( \pi_{f} \right ) -\sum_{i=1}^{n} p_{i,j}\cdot \sum_{k=1}^{n} y_{k,y} \right ]\cdot SP \cdot \varepsilon_{elec} 
	\setlength	\belowdisplayskip {3pt} 
\end{equation}

The carbon emissions during machine idle time equal the product of the total idle time of all machines, the rated energy consumption during idling, and the carbon emission coefficient for electricity.

(3) Carbon emissions from auxiliary materials consumption ($CE_{au}$):
\begin{equation}
	\setlength	\abovedisplayskip {3pt} 
	\begin{split}
		CE_{au} &=\sum_{j=1}^{n} \left ( \frac{p_{i,j} }{T_{i,j}^{cool} } \cdot IC_{i,j}^{cool}  \cdot EF_{i,j}^{cool} +\frac{p_{i,j} }{T_{i,j}^{lu} } \cdot LO_{i,j}^{lu}  \cdot EF_{i,j}^{lu}\right ) \\
		&= \sum_{j=1}^{n_{i}}p_{i,j} \cdot\left ( \frac{IC_{i,j}^{cool}\cdot EF_{i,j}^{cool}}{T_{i,j}^{cool}} + \frac{LO_{i,j}^{lu}\cdot EF_{i,j}^{lu}}{T_{i,j}^{lu}}\right ) \\
		&=\sum_{j=1}^{n_{i}}p_{i,j}\cdot \varepsilon_{elec}
	\end{split}
	\setlength	\belowdisplayskip {3pt} 
\end{equation}

Machines in the manufacturing process consume auxiliary materials, particularly coolants and lubricating oils. Coolants are typically circulated by a cooling pump and gradually diminish due to adherence to metal chips until replenished. Lubricating oil is primarily used for machine tool spindles and guide rails, where a small amount of oil is injected into the spindle and guide rails at regular intervals. To simplify the model, the time for tool changes and the recycling process for lubricating oil and coolant are ignored. \textcolor[rgb]{0,0,0}{In our work,  $T_{i,j}^{cool}$ represents the average interval for refreshing the coolant in the machine processing job $O_{i,j}$}. $EF_{i,j}^{cool}$  is the carbon emission coefficient for the machine’s coolant. $p_{i,j}$  indicates the processing time of job  $i$ on machine $j$. $T_{i,j}^{lu}$  denotes the average interval for discharging lubricating oil in the machine processing job $O_{i,j}$.  $LU_{i,j}^{lu}$ represents the initial amount of lubricating oil for the machine processing job $O_{i,j}$. $EF_{i,j}^{lu}$  is the carbon emission coefficient for the machine’s lubricating oil $IC_{i,j}^{cool}$  indicates the initial amount of coolant for the machine processing job $O_{i,j}$. $\varepsilon_{elec}$  is the carbon emission coefficient for auxiliary materials.

\textbf{Decision variables:}
\begin{equation}
	\setlength	\abovedisplayskip {3pt} 
	x_{i,k} =\left\{\begin{matrix}
		1 & ,When \ job \ i \ is \ at \ position \ k\\
		0&, else
	\end{matrix}\right.  
	\setlength	\belowdisplayskip {3pt} 
\end{equation}
\begin{equation}
	y_{k,f} =\left\{\begin{matrix}
		1 & ,The \ job \ at \ position \  k \\
		& is \ processed \  in \ factory \  f\\
		0&, else
	\end{matrix}\right.  
\end{equation}

\textbf{Constraints:}
\begin{equation}
	\sum_{k=1}^{n}x_{i,k} =1,i\in \left \{ 1,2,...,n \right \} 
\end{equation}
\begin{equation}
	\setlength	\abovedisplayskip {3pt} 
	\sum_{i=1}^{n}x_{i,k} =1,k\in \left \{ 1,2,...,n \right \} 
	\setlength	\belowdisplayskip {3pt} 
\end{equation}
\begin{equation}
	\setlength	\abovedisplayskip {3pt} 
	\sum_{f=1}^{F}y_{k,f} =1,k\in \left \{ 1,2,...,n \right \} 
	\setlength	\belowdisplayskip {3pt} 
\end{equation}
\begin{equation}
	\setlength	\abovedisplayskip {3pt} 
	C_{k,f}\ge C_{k,j-1}+  \sum_{i=1}^{n}x_{i,k}\cdot p_{i,j},k\in \left \{ 1,2,...,n \right \} ,j\in \left \{ 2,3,...,m \right \} 
	\setlength	\belowdisplayskip {3pt} 
\end{equation}
\begin{equation}
	\setlength	\abovedisplayskip {3pt} 
	\begin{split}
		C_{k,j}\ge C_{l,j}+  \sum_{i=1}^{n}x_{i,k}\cdot p_{i,j}- M\left (  1-y_{k,f}\right ) - M\left (  1-y_{l,f}\right ),\\
		k\in \left \{ 2,3,...,m \right \}, l> k,i\in \left \{ 1,2,...,m \right \},f\in \left \{ 1,2,...,F \right \}  
	\end{split} 
	\setlength	\belowdisplayskip {3pt}	
\end{equation}
\begin{equation}
	\setlength	\abovedisplayskip {3pt} 
	C_{max}\ge C_{k,m},k\in \left \{ 1,2,...,n \right \}  
	\setlength	\belowdisplayskip {3pt} 
\end{equation}
\begin{equation}
	\setlength	\abovedisplayskip {3pt} 
	C_{k,j}\ge 0,k\in \left \{ 1,2,...,n \right \},j\in \left \{ 1,2,...,m \right \}  
	\setlength	\belowdisplayskip {3pt} 
\end{equation}

\textcolor[rgb]{0,0,0}{Equation (1) is the objective function, minimizing the makespan and total carbon emissions. Equations (2) and (3) define the range of values for the decision variables. Equations (4) and (5) ensure that each job appears only once. Equation (6) ensures that each job is assigned to only one factory. Equation (7) maintains the constraint of adjacent operations for the same job, ensuring a job is not processed on multiple machines simultaneously. Equation (8) enforces the constraint for adjacent jobs on the same machine, preventing simultaneous processing of multiple jobs on one machine. Equation (9) indicates that the makespan must be greater than, or equal to, the completion time of the job on the last machine in all factories. Equation (10) states that the completion time for all operations must be greater than zero.}

\section{The Knowledge-driven Memetic Algorithm}
\label{KDMA}
MA is an algorithmic framework inspired by Darwin’s principle of natural evolution and Dawkin’s concept of memes, combining global population-based search with individual heuristic local search \cite{wu2022ensemble}. MAs utilize a mechanism that combines global and local searches, making their search efficiency several orders of magnitude faster than traditional genetic algorithms in certain problem domains. Different MAs can be constructed using various global search strategies and local search strategies \cite{song2023learning,LOU2022101204}.

\begin{figure}[htp]
\centering
\includegraphics[width=0.45\textwidth]{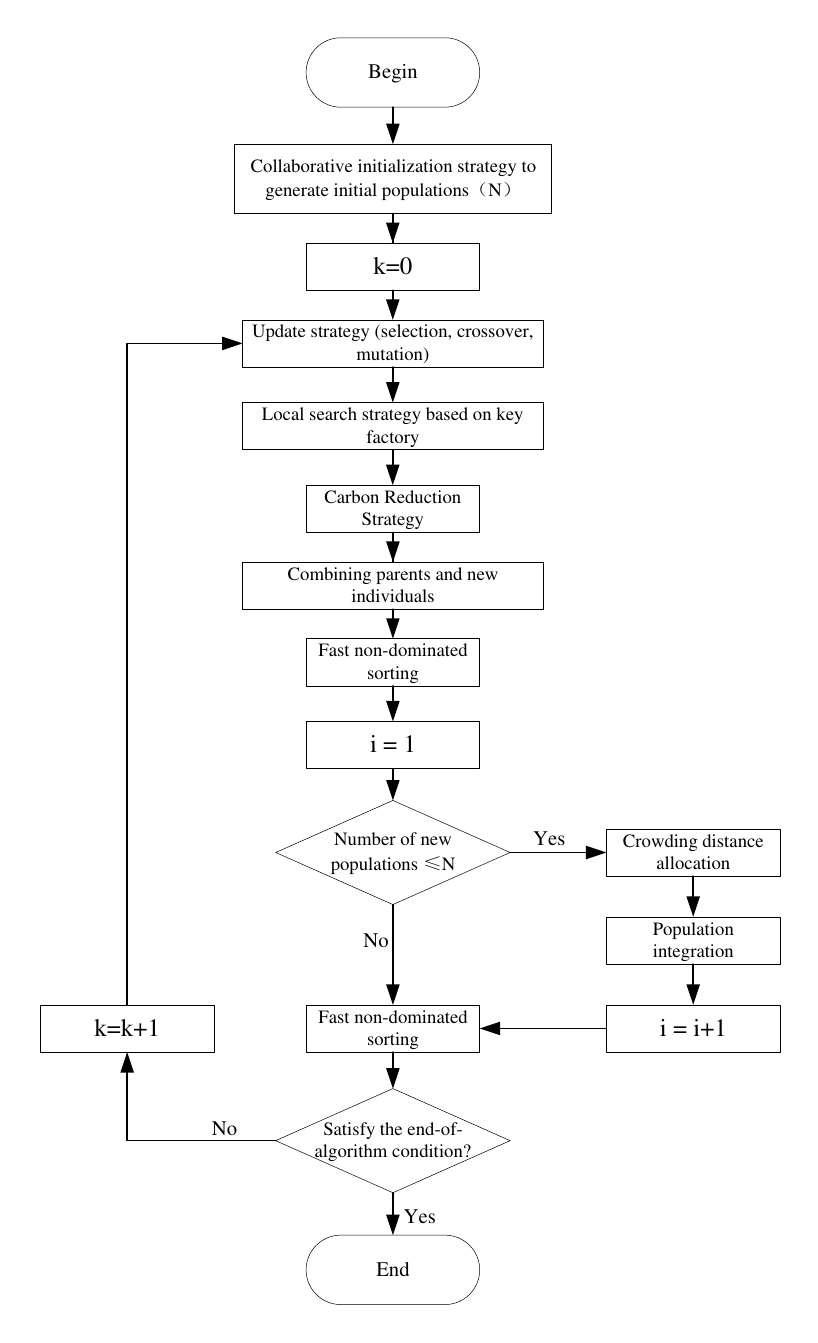}
\caption{Flowchart through the KDMA}
\label{Flowchart through the KDMA}
\end{figure}

\textcolor[rgb]{0,0,0}{Knowledge-driven strategies are an effective way to improve on MA. These strategies can use features in the problem that drive the algorithm search. In our study, a knowledge-driven MA (KDMA) is proposed to solve the EEDHFSSP. This algorithm primarily consists of a collaborative initialization strategy, an updating strategy, a local search strategy based on key factories, and a carbon reduction strategy. The algorithm flowchart is illustrated in Figure.\ref{Flowchart through the KDMA}.}

\subsection{Collaborative Initialization Strategy}

The quality of initial solutions in combinatorial optimization problems significantly affects the algorithm’s performance \cite{song2024generalized}. Optimizing initial solutions can improve algorithm performance \cite{wu2022coordinated}. Generally, it is challenging to find an optimal solution that satisfies dual objectives, namely, minimizing makespan and total carbon emissions. To ensure a more uniform distribution of initial solutions, we propose a collaborative initialization strategy based on the encoding characteristics of EEDHFSSP. The Nawaz-Enscore-Ham (NEH) heuristic algorithm has been proven to be one of the most effective heuristic algorithms for addressing flow shop scheduling problems \cite{Fernandez2017}. A modified NEH (MNEH) algorithm is introduced to generate an initial solution with the minimum makespan during initialization. The steps in the MNEH algorithm are as follows:

\textbf{Step 1:} Generate a random initial job allocation sequence $\Pi=\left [ O_{1},...,...,O_{n} \right ] $.

\textbf{Step 2:} Sequentially insert each job $O_{i}\left ( i=1,2,...,n\right )$ into all possible positions in each factory and select the scheduling sequence with the minimum makespan value.

\textbf{Step 3:} The scheduling sequence after inserting job  $O_{i}$ is used for the insertion of the next job, continuing until all jobs are allocated to factories.

To minimize total carbon emissions, the following strategy is adopted to generate an initial solution with the minimum total carbon emissions:

\textbf{Step 1:} Sort all jobs in non-ascending order of their rated processing power, denoted by  $ O =\left [ O_{1},...,...,O_{n} \right ]$.

\textbf{Step 2:} Following this order, sequentially insert each job $O_{i}\left ( i=1,2,...,n \right ) $  into all possible positions and calculate the carbon emissions associated with each position.

\textbf{Step 3:} Select the position with the minimum carbon emissions as the final allocation for each job. Continue this process until all jobs are allocated. The remaining initial solutions are generated randomly to maintain diversity in the population.

\textcolor[rgb]{0,0,0}{Inevitably, there will be situations where constraints are difficult to satisfy during the decoding process. The algorithm attempts to backtrack the operation in case of a constraint violation, until the constraint is met.}

\subsection{Variation Strategy}

During the algorithm’s search process, crossover and mutation operators are used to update solutions. Considering the characteristics of EEDHFSSP, a tournament selection method is adopted to choose parents. The partially mapped crossover (PMX) operator, widely applied in scheduling domains, is used for crossover operations, and a swap mutation operator is employed for mutation \cite{song2023rl}.

Tournament selection is a local competition-based selection method. From the population, $k$  individuals are randomly selected for comparison, and the individual with the highest fitness value is chosen to enter the parent population. This process is repeated $N$  times until the parent population reaches the desired size.

\begin{figure}[htp]
\centering
\includegraphics[width=0.45\textwidth]{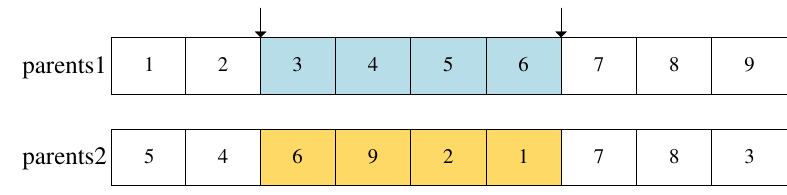}
\caption{PMX operator selection operation}
\label{PMX operator selection operation}
\end{figure}

The steps for the PMX operator are as follows:

\textbf{Step 1:} Randomly select two positions on the chosen parent chromosomes, and define the elements between these two positions as the matched sub-string, as shown in Figure \ref{PMX operator selection operation}.

\begin{figure}[htp]
\centering
\includegraphics[width=0.45\textwidth]{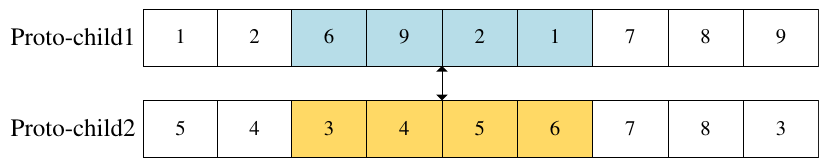}
\caption{PMX operator swap operation}
\label{PMX operator swap operation}
\end{figure}

\textbf{Step 2:} Swap the two matched sub-strings of the parent chromosomes to obtain two temporary offspring chromosomes, as shown in Figure \ref{PMX operator swap operation}.

\begin{figure}[htp]
\centering
\includegraphics[width=0.45\textwidth]{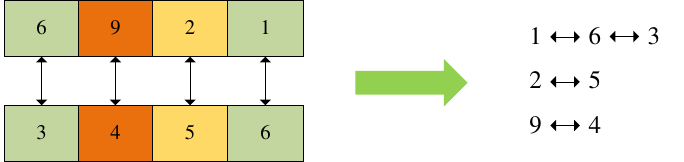}
\caption{PMX operator determining conflicting jobs operation}
\label{PMX operator determining conflicting jobs operation}
\end{figure}

\textbf{Step 3:} Determine and map the relationships of conflicting jobs based on the swapped pairs of genes, as shown in Figure \ref{PMX operator determining conflicting jobs operation}.

\begin{figure}[htp]
\centering
\includegraphics[width=0.45\textwidth]{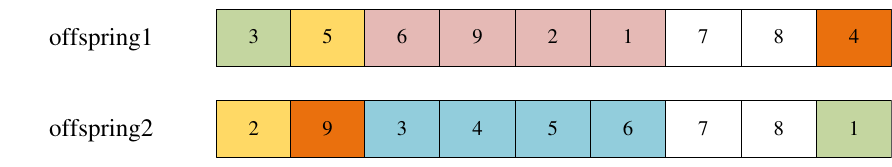}
\caption{PMX operator determining job order operation}
\label{PMX operator determining job order operation}
\end{figure}

\textbf{Step 4:} Ensure the feasibility of the job sequence based on the mapping relationship, with no need for any changes to the sub-string, as shown in Figure \ref{PMX operator determining job order operation}.

The swap mutation operation involves randomly selecting two genes on the chosen parent chromosome and swapping the jobs on these two genes.

\subsection{Local Search Strategy Based on Key Factories}

According to existing research, local search heuristics are one of the effective methods for solving single-objective optimization problems (SOPs) \cite{song2024energy,zhao2022reinforcementa}. However, in multi-objective scheduling problems, conflicts exist among multiple objectives, and optimizing one objective may lead to the deterioration of another. Therefore, local search heuristics effective for single-objective optimization problems cannot be directly applied to EEDHFSSP. By analyzing and summarizing the characteristics of EEDHFSSP, a local search strategy based on key factories is designed to enhance the performance of the algorithm.

Four types of local search operation operators are designed. These four local search operators can be divided into operations within factories and operations between factories. The operations within factories include $Insert\_in\_Crifac$  and $Swap\_in\_Crifac$, while operations between factories include $Insert\_between\_fac$  and $Swap\_between\_fac$ , as follows:

(1) $L1\left ( Insert\_in\_Crifac \right ) $: job insertion within key factories. In a key factory, randomly select a job and insert it into all possible positions within that key factory.

(2)  $L2\left ( Swap\_in\_Crifac \right ) $: job swap within key factories. In a key factory, randomly select two positions $j$  and $k$ , and swap the jobs at positions $j$  and $k$.

(3) $L3\left ( Insert\_between\_fac \right ) $: job exchange between factories. Randomly select a non-key factory $b$, choose a position $k$  within factory $b$, then select a position  $j$ within a key factory, and swap the jobs at positions $k$  and $j$.

(4) $L4\left ( Swap\_between\_fac \right ) $: job insertion between factories. In a key factory, randomly select a position $j$  and remove the job from that position, then insert it into any position $k$ in any non-key factory $b$.

\subsection{Carbon Reduction Strategy}

In the DHFSSP, the total carbon emissions consist of three parts: carbon emissions produced during job processing, emissions generated during machine idle time, and emissions arising from consumption of auxiliary materials. Since the processing time for each job is the same in each factory, the carbon reduction strategy designed here mainly targets the carbon emissions generated during machine idle time.

\begin{figure}[htp]
\centering
\includegraphics[width=0.5\textwidth]{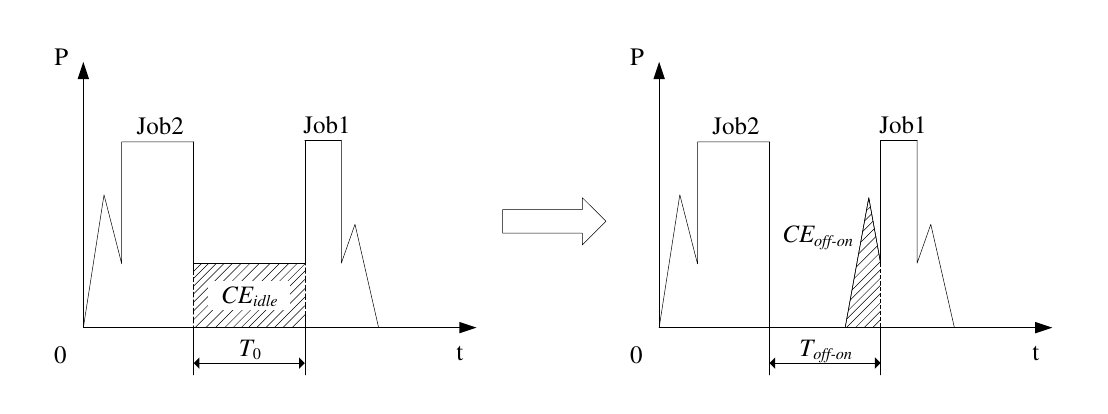}
\caption{Example of the machine on/off strategy}
\label{Example of the machine on/off strategy}
\end{figure}

Assume  $T_{0}$ is the idle waiting time of the machine, $T_{off-on}$  is the time required to turn off and then turn on the machine, $CE_{idle}$ denotes the carbon emissions produced during machine idle time, and $CE_{off-on}$  is the carbon emissions generated by turning off and then on the machine. When $CE_{idle} > CE_{off-on}$  and $T_{0} > T_{off-on}$, shutting down the idle machine may reduce the carbon emissions generated during idle time, as illustrated in Figure \ref{Example of the machine on/off strategy}.

\section{Experimental Results and Analysis}
\label{Experiment}
\subsection{Experimental Setup}

To validate the effectiveness of KDMA in solving the EEDHFSSP, a series of comparative experiments were designed and conducted on a PC with Microsoft Windows 10 operating system, equipped with an Intel Core I5-10500 CPU @ 3.10 GHz and 16 GB RAM.

As the EEDHFSSP is a new problem without open datasets, a set of test datasets were established according to the characteristics and main parameters of EEDHFSSP. The instances include different numbers of factories $f=\left \{ 2,3,4,5,6 \right \} $, jobs $n=\left \{ 20,50,100 \right \} $, and machines $m=\left \{ 2,5,8 \right \} $. For each $\left ( f,n,m \right )$  combination, 10 instances were generated, giving a total of 450 instances. For each instance, the standard processing time of jobs $p_{i,j}$ was uniformly sampled from the discrete range $\left [ 10,50 \right ]$, the unit energy consumption of jobs during processing on machines  $PP_{i,j}$ was sampled from $\left [ 5,10 \right ]$, the idle power consumption of each machine $SP=2$, and the carbon emissions for turning off/on a machine were set to $CE_{off-on}=6$. Following the Corporate Greenhouse Gas Emission Accounting and Reporting Guide-Power Facilities (Revised Edition 2022), the electricity carbon emission coefficient  was set. Additionally, the carbon emission coefficient \textcolor[rgb]{0,0,0}{ $ \varepsilon _{elec}=0.581 kgCO_{2}e/kWh $}for auxiliary materials used during machine $ \varepsilon _{au}$ processing was uniformly sampled from\textcolor[rgb]{0,0,0}{ $\left [ 0.05,0.1 \right ]\left ( kgCO_{2}e/kg \right) $ }, as it is related to the type of machine used.

\textcolor[rgb]{0,0,0}{To evaluate the performance of the algorithm, the Pareto solution set from various aspects, including diversity of solutions, convergence, and overall performance, using three different performance metrics:}

(1) $Spread$: This metric evaluates the extent of distribution of the Pareto solution set obtained by the algorithm in the objective space.

(2) $GD$ (Generation Distance): This metric calculates the average minimum distance of each point in the solution set $P$ to the reference set $P^{\ast } $, reflecting the convergence of the obtained solution set. A smaller $GD$ value indicates better convergence of the algorithm.
\begin{equation}
\setlength \abovedisplayskip {3pt} 
GD\left ( P,P^{\ast }  \right ) =\frac{\sqrt{\sum_{y\in P^{\ast } } \min_{x\in P}  dis \left ( x,y \right ) ^{2}  }  }{\left | P \right | }  
\setlength	\belowdisplayskip {3pt} 
\end{equation}

where, $P$ is the Pareto solution set obtained by the algorithm, $P^{\ast}$ denotes a set of uniformly distributed reference points sampled from the Pareto front, and $dis \left ( x,y \right ) $ represents the Euclidean distance between point $x$ in the solution set $P$  and point $P^{\ast }$ in the reference set.

(3) $IGD$ (Inverted Generation Distance): This metric calculates the average distance from each reference point to the nearest solution, reflecting the overall performance of the algorithm. A smaller $IGD$ value indicates better overall performance.
\begin{equation}
\setlength	\abovedisplayskip {3pt} 
IGD\left ( P,P^{\ast }  \right ) =\frac{\sum_{y\in P } \min_{x\in P^{\ast }}  dis \left ( x,y \right )}{\left | P^{\ast } \right | }   
\setlength	\belowdisplayskip {3pt} 
\end{equation}

where, $P$ is the Pareto solution set obtained by the algorithm, $P^{\ast}$  is a set of uniformly distributed reference points sampled from the Pareto front, and $dis \left ( x,y \right )$ represents the Euclidean distance between point $x$ in the reference set $P^{\ast}$ and point $y$ in the solution set $P$.

\subsection{Parameter Settings}

The KDMA involves three key parameters: (1) population size ($PS$), (2) crossover probability ($p_{c}$), and (3) mutation probability ($p_{m}$). To study the effects of these parameters on the performance of the KDMA, this research used the widely utilized Taguchi experimental design method. As shown in Table \ref{Parameter level correspondence}, four different levels were set for each parameter: $PS=\left \{ 50,100,150,200 \right \} $ , $p_{c}=\left \{ 0.7,0.8,0.9,1.0 \right \} $ , and $p_{m}=\left \{ 0.1,0.2,0.3,0.4 \right \} $ . The orthogonal array $L_{16}\left ( 4^{3}  \right )$ generated from these settings was used in the parameter calibration experiments. Additionally, the overall performance metric was $IGD$ adopted to measure the effects of different parameter configurations on the dataset of instances. For each instance, the KDMA was executed 10 times for each parameter configuration, with the maximum number of function evaluations set to 25000.

\begin{table}[!ht]
	\centering
	\caption{Parameter level correspondence}
	\label{Parameter level correspondence}
	\footnotesize
 \setlength{\tabcolsep}{0.5cm}
	\begin{tabular}{|c|c|c|c|c|} \hline
		\multicolumn{1}{|c|}{\multirow{2}{*}{Parameter}} & \multicolumn{4}{c|}{Level}\\ \cline{2-5}
		& 1 & 2 & 3 & 4\\ \hline
	$PS$ &50 &100 &150 &200\\ \hline
	$p_{c}$  &0.7 &0.8 & 0.9&1.0 \\\hline
		$p_{m}$ & 0.1 & 0.2 & 0.3& 0.4 \\ \hline
	\end{tabular}
\end{table}

\begin{figure}[htp]
\centering
\includegraphics[width=0.45\textwidth]{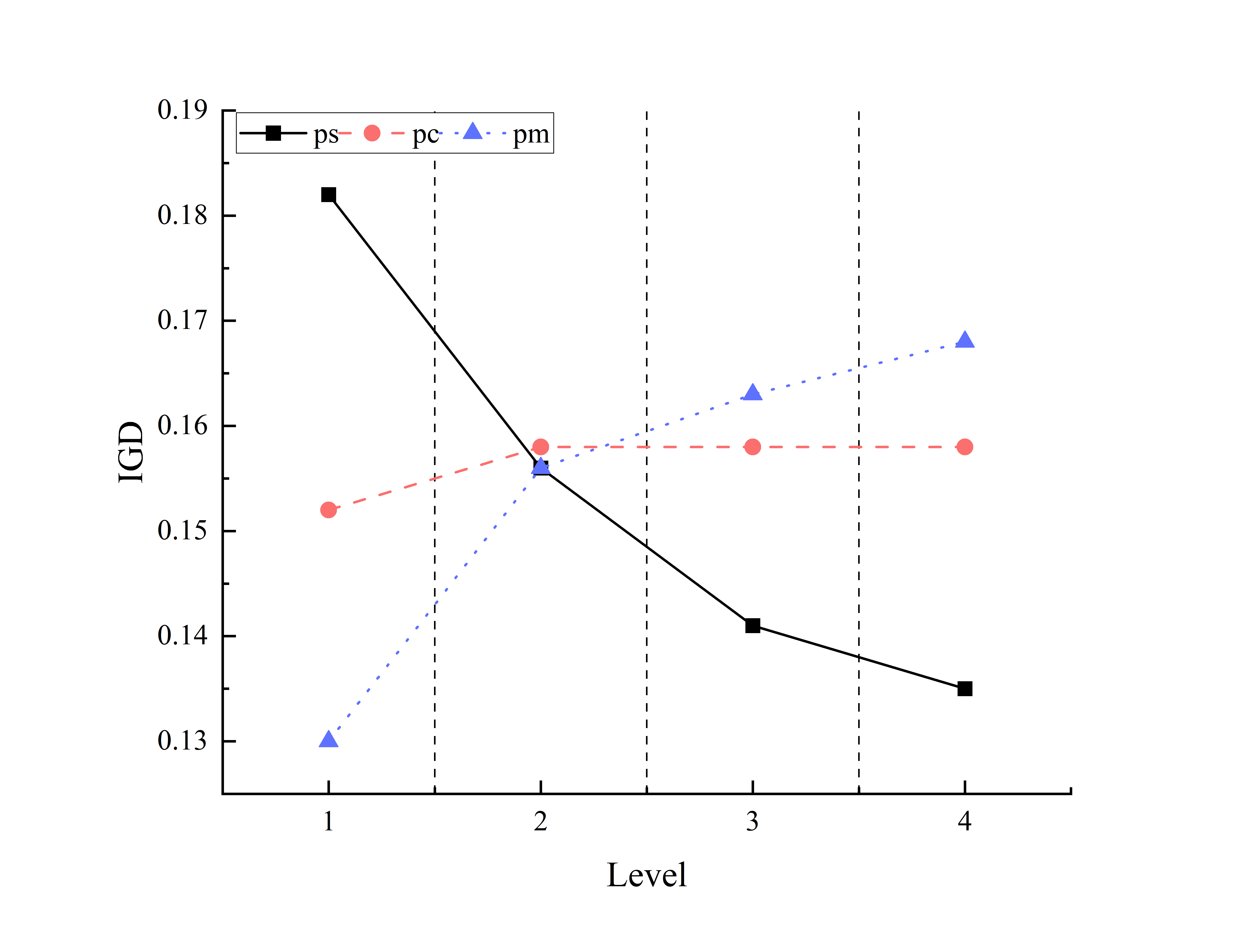}
\caption{Comparison of parameter performance metrics at different levels}
\label{Comparison of parameter performance metrics at different levels}
\end{figure}

Table \ref{Average metric values for each parameter} presents the significance levels of the metrics and parameters, where the   value indicates the maximum difference between the average metric value and other metric values. A larger $delta$ value implies that the corresponding parameter is more important. Figure \ref{Comparison of parameter performance metrics at different levels} shows the performance metrics of the three parameters at different levels.

\begin{table}[htbp]
	\centering
	\caption{Average metric values for each parameter}
	\label{Average metric values for each parameter}
	\begin{tabularx}{0.45\textwidth}{m{1.5cm}m{1.5cm}m{1.5cm}m{1.5cm}}
		\toprule[1.5pt]
		Level & $PS$ & $p_{c}$ & $p_{m}$\\
		\midrule[0.75pt]
		1 & 0.181088 & 0.152377 & 0.131506 \\
		2 & 0.157780 & 0.155338 & 0.155154 \\
		3 & 0.142486 & 0.154651 & 0.164213 \\
            4 & 0.135692 & 0.154680 & 0.166173 \\
            $Average$ & 0.154262 & 0.154262 & 0.154262 \\
            $delta$ & 0.045396 & 0.002961 & 0.034667 \\
            Rank & 1 & 3 & 2 \\
		 \bottomrule[1.5pt]                        
	\end{tabularx}
\end{table}

As presented in Table \ref{Average metric values for each parameter}, population size ($PS$) has the greatest effect on the algorithm’s performance, followed by mutation probability ($p_{m}$), and crossover probability ($p_{c}$) has the least impact. This means that the importance ranking is: population size ($PS$) as the most important, followed by mutation probability ($p_{m}$), and then crossover probability ($p_{c}$). According to Figure \ref{Comparison of parameter performance metrics at different levels} , the optimal values for the three key parameters are: $PS=100,p_{c}=0.9,p_{m}=0.2$.

\subsection{Experiment on the Effectiveness of Strategies}

\begin{figure}[htp]
\centering
\includegraphics[width=0.45\textwidth]{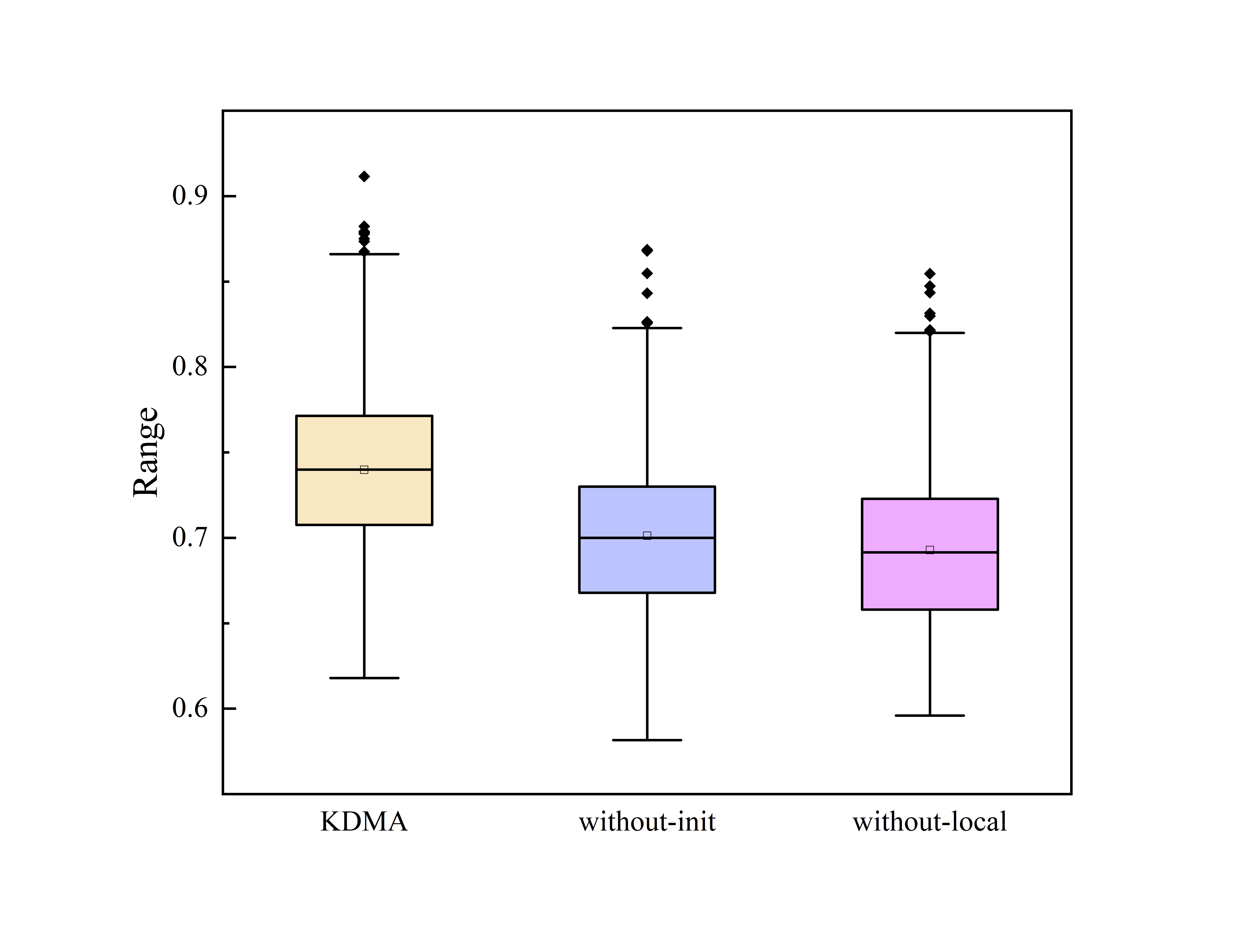}
\caption{Comparison of $SP$ values in the strategy effectiveness experiments}
\label{Comparison of $SP$ values in the strategy effectiveness experiment}
\end{figure}

\begin{figure}[htp]
\centering
\includegraphics[width=0.45\textwidth]{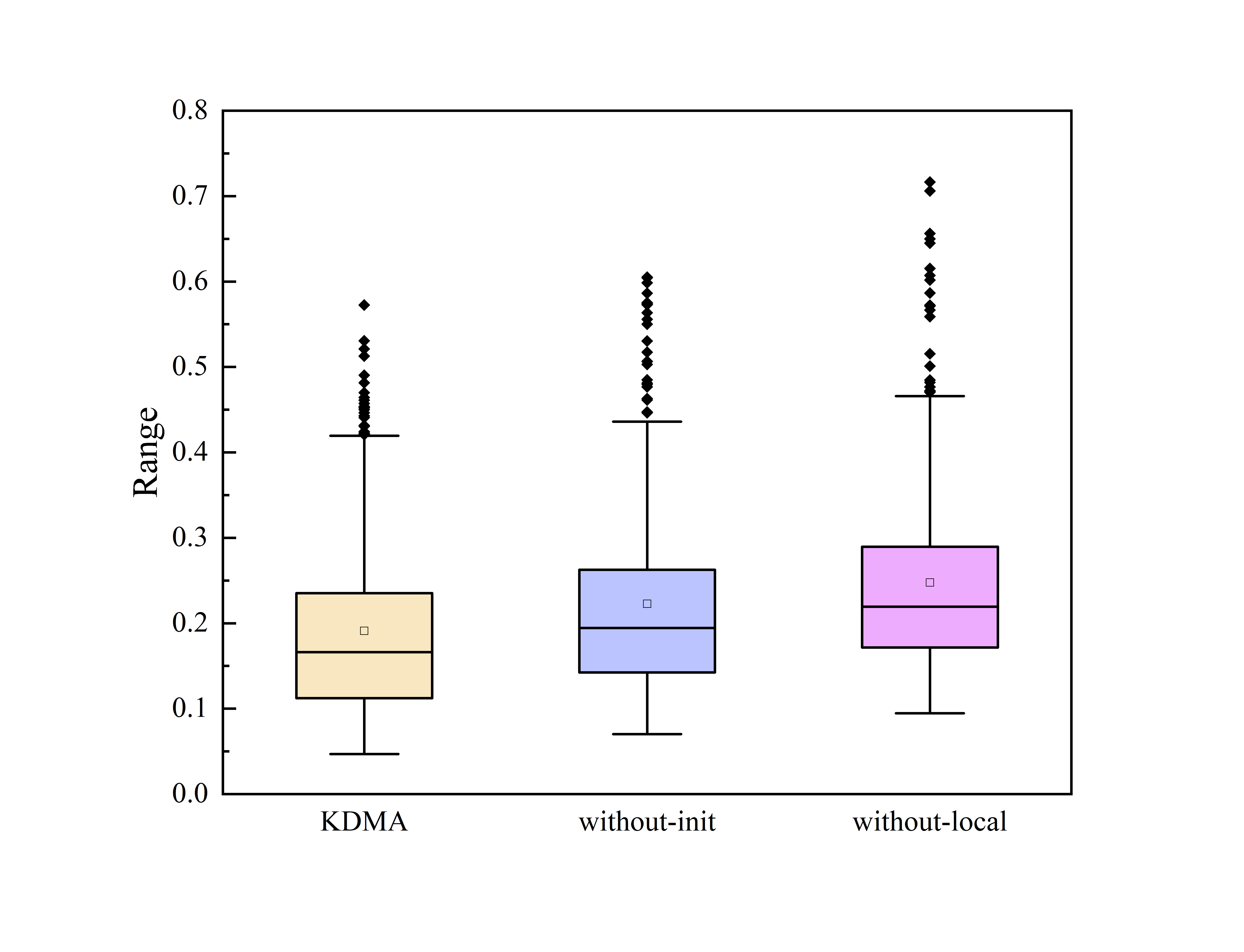}
\caption{Comparison of $GD$ values in the strategy effectiveness experiment}
\label{Comparison of $GD$ values in the strategy effectiveness experiment}
\end{figure}

\begin{figure}[htp]
\centering
\includegraphics[width=0.45\textwidth]{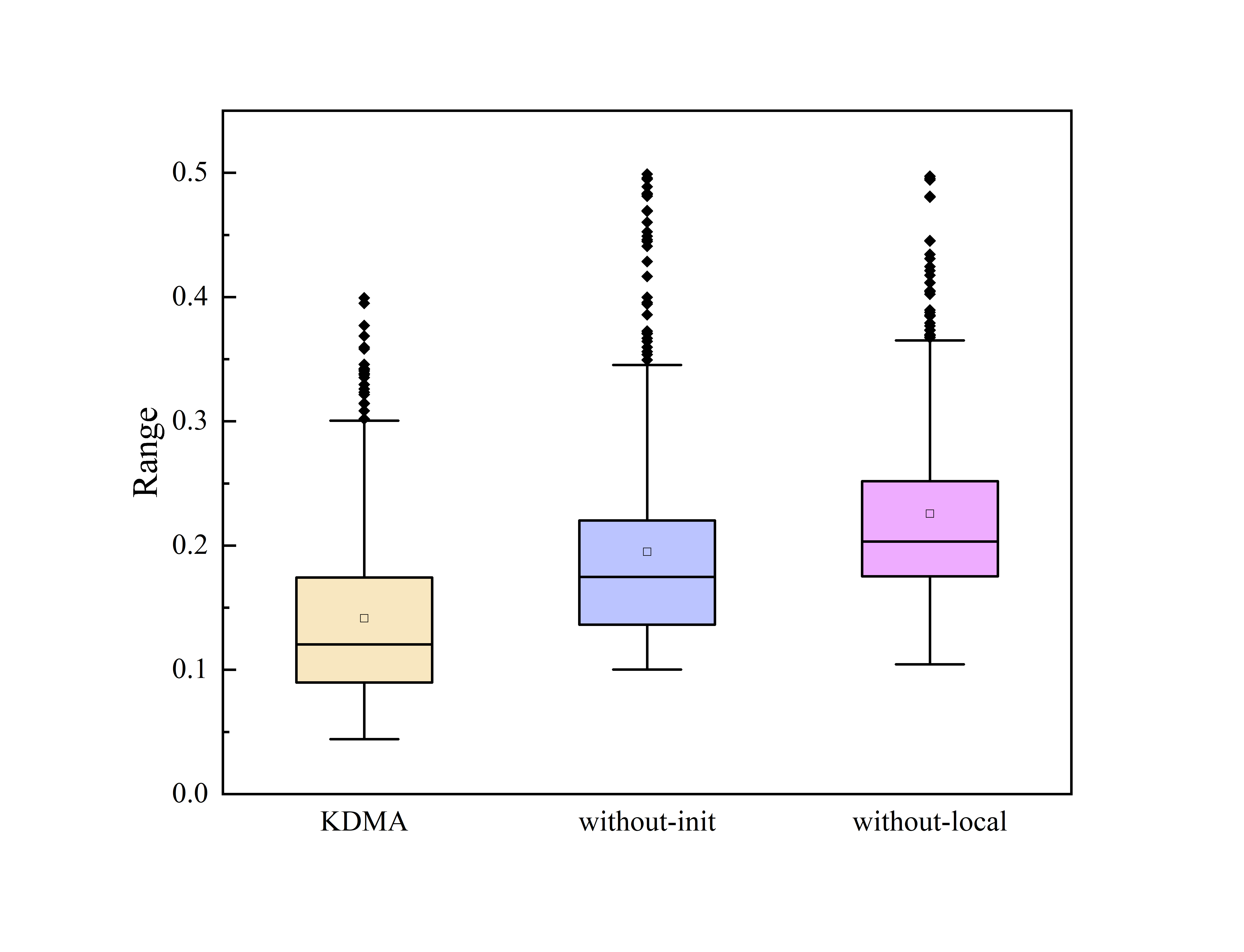}
\caption{Comparison of $IGD$ values in the strategy effectiveness experiment}
\label{Comparison of $IGD$ values in the strategy effectiveness experiment}
\end{figure}

To validate the efficacy of the proposed collaborative initialization strategy and the local search strategy based on key factories, this study conducted comparative experiments among KDMA with both the collaborative initialization and local search strategies, KDMA without the collaborative initialization strategy (wherein the initial population is entirely generated randomly, denoted as “\textbf{without-init}”), and KDMA without the knowledge-based local search strategy (referred to as “\textbf{without-local}”). Each algorithm was independently run 10 times on each instance, all using the same termination criterion (the maximum number of function evaluations set to 25,000). Figures \ref{Comparison of $SP$ values in the strategy effectiveness experiment},\ref{Comparison of $GD$ values in the strategy effectiveness experiment},\ref{Comparison of $IGD$ values in the strategy effectiveness experiment} present the $Spread$, $GD$, and $IGD$ values obtained after running the three algorithms: KDMA, without-init, and without-local.

As illustrated in Figure \ref{Comparison of $SP$ values in the strategy effectiveness experiment}, across all instances, KDMA’s results in the $SP$ metric significantly surpass those of without-init, while the difference with without-local is less pronounced. This indicates that the collaborative initialization strategy proposed herein contributes to enhancing the diversity of solutions. As shown in Figure \ref{Comparison of $GD$ values in the strategy effectiveness experiment}, KDMA outperforms the two strategies-lacking algorithms in the $GD$ performance metric, suggesting that both the collaborative initialization strategy and the knowledge-based local search strategy improve the convergence of the algorithm. As depicted in Figure \ref{Comparison of $IGD$ values in the strategy effectiveness experiment}, regarding the $IGD$ performance metric, KDMA’s results are markedly superior to both without-init and without-local, demonstrating that the strategies proposed herein enhance the overall performance of the algorithm.

\begin{figure}[htp]
\centering
\includegraphics[width=0.45\textwidth]{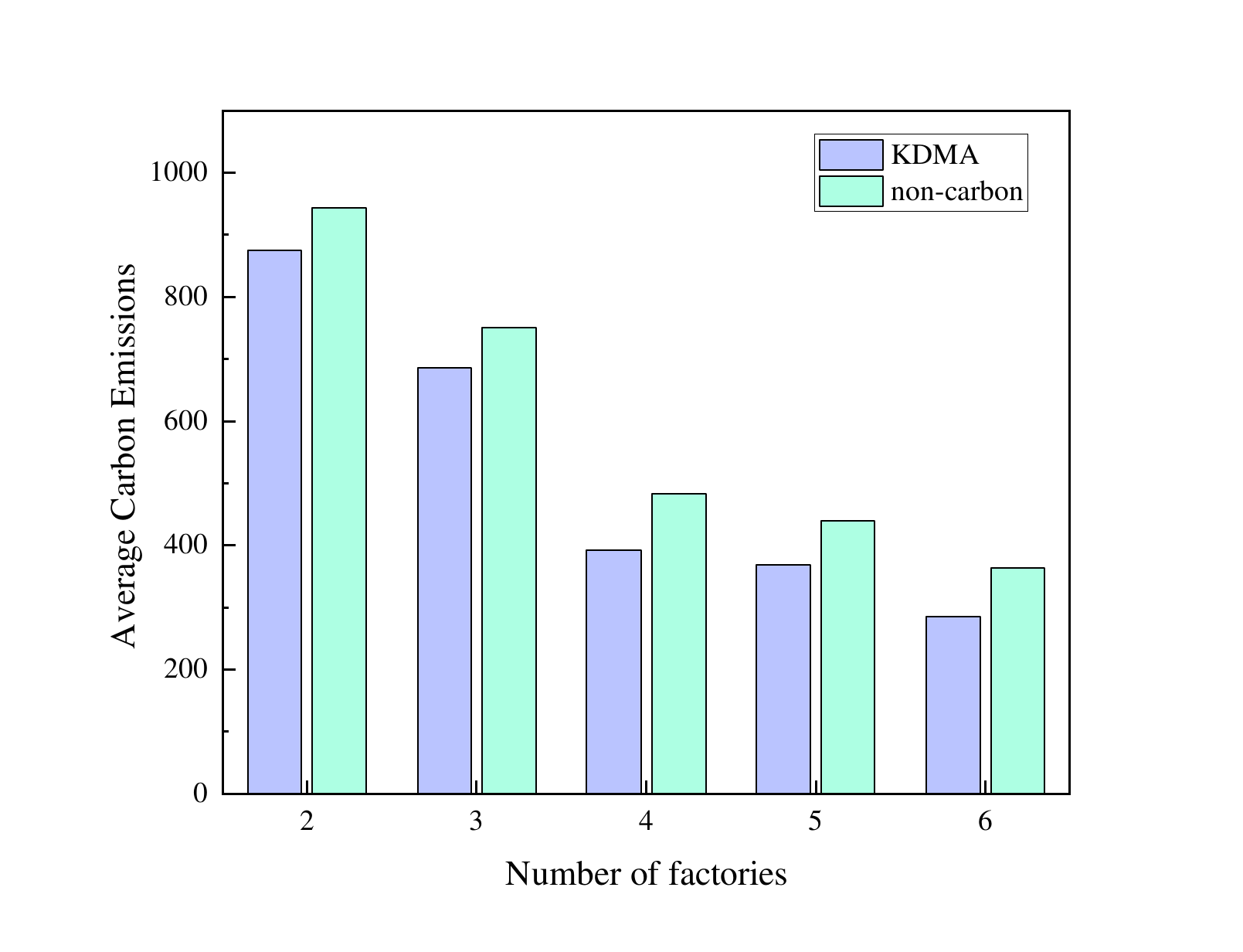}
\caption{\textcolor[rgb]{0,0,0}{Comparison of average carbon emissions for different numbers of factories}}
\label{Comparison of average carbon emissions for different numbers of factories}
\end{figure}

To validate the effectiveness of the carbon reduction strategy in KDMA, KDMA was compared with the carbon reduction strategy to KDMA without the carbon reduction strategy (referred to as non-carbon). Each algorithm was independently run 10 times on each instance, all using the same termination criterion (the maximum number of function evaluations set to 25,000). The total carbon emission results for both algorithms were compiled, and the average carbon emissions for each type of factory were calculated based on different numbers of factories. Figure \ref{Comparison of average carbon emissions for different numbers of factories} shows the average carbon emissions for the two algorithms across different numbers of factories.

As shown in Figure \ref{Comparison of average carbon emissions for different numbers of factories}, regardless of the number of factories, the average carbon emission values for KDMA with the carbon reduction strategy are consistently lower than those for non-carbon schema, which do not use the carbon reduction strategy. This indicates that the proposed carbon reduction strategy reduces the carbon emissions in distributed homogeneous flow shop scheduling problems.

\subsection{Comparison of Different Algorithms}

Three classic multi-objective optimization algorithms were selected: a non-dominated sorting genetic algorithm-II (NSGA-II) \cite{996017}, a multi-objective evolutionary algorithm based on decomposition (MOEA/D) \cite{6252954}, and a weighting achievement scalarizing function genetic algorithm (WASFGA) \cite{Ruiz2015}, for comparison with KDMA. The experiments were conducted on a generated instance dataset, with identical algorithm parameters and termination conditions set for all. Each of the four algorithms was independently run 10 times on each  $\left ( f,n,m \right )$ combination dataset, and the results were averaged over these 10 experiments. The efficacy of KDMA in solving the EEDHFSSP was verified by comparing the $SP$ , $GD$ , and  $IGD$ metrics. Tables \ref{Comparison of average $SP$ values}-\ref{Comparison of average $IGD$ values} respectively list the  $SP$ , $GD$ , and  $IGD$ values for KDMA, NSGA-II, MOEA/D, and WASFGA.

\begin{table*}[!ht]
	\center
	\caption{Comparison of average $SP$  values}
	\label{Comparison of average $SP$  values}
		\setlength{\tabcolsep}{0.8cm}
			\begin{tabular}{|c|c|c|c|c|c|c|} \hline
			$F$ & $m$ & $n$ & KDMA & NSGA-II & MOEA/D & WASFGA  \\ \hline
				
				\multirow{9}{*}{2}&\multirow{3}{*}{2} & 20 &	\textbf{0.803659}	&0.743158	&0.752892&0.653639  \\ \cline{3-7}
				\multirow{9}{*}{}&\multirow{3}{*}{} & 50	&\textbf{0.762419} &0.715467&	0.748665 &	0.643803 \\ \cline{3-7}
				\multirow{9}{*}{}&\multirow{3}{*}{} & 100	&\textbf{0.782145} &0.718958&	0.755771 &	0.669015 \\ \cline{2-7}
				
				\multirow{9}{*}{ }&\multirow{3}{*}{5} & 20 &	\textbf{0.774332}	&0.688873	&0.767472&0.600316  \\ \cline{3-7}
				\multirow{9}{*}{}&\multirow{3}{*}{} & 50	&\textbf{0.751021} &0.695362&	0.715166 &	0.643988 \\ \cline{3-7}
				\multirow{9}{*}{}&\multirow{3}{*}{} & 100	&\textbf{0.780631} &0.702863&	0.753433 &	0.7090495 \\ \cline{2-7}
				
				\multirow{9}{*}{ }&\multirow{3}{*}{8} & 20 &	\textbf{0.767921}	&0.671917	&0.698723&0.671018  \\ \cline{3-7}
				\multirow{9}{*}{}&\multirow{3}{*}{} & 50	&\textbf{0.747668} &0.696998&	0.706409 &	0.636901 \\ \cline{3-7}
				\multirow{9}{*}{}&\multirow{3}{*}{} & 100	&\textbf{0.753020} &0.695810&	0.724216 &	0.669222 \\ \hline
				
				\multirow{9}{*}{3}&\multirow{3}{*}{2} & 20 &	\textbf{0.723902}	&0.705576	&0.719963&0.632921  \\ \cline{3-7}
				\multirow{9}{*}{}&\multirow{3}{*}{} & 50	&\textbf{0.763501} &	0.712029&	0.742161&	0.605642 \\ \cline{3-7}
				\multirow{9}{*}{}&\multirow{3}{*}{} & 100	&\textbf{0.738157} &	0.716682&	0.717895&	0.618383 \\ \cline{2-7}
				
				\multirow{9}{*}{ }&\multirow{3}{*}{5} & 20 &	\textbf{0.731253}	 &	0.688978 &	0.691626 &	0.618832  \\ \cline{3-7}
				\multirow{9}{*}{}&\multirow{3}{*}{} & 50	&\textbf{0.747066	} &0.696465&	0.718189&	0.600427 \\ \cline{3-7}
				\multirow{9}{*}{}&\multirow{3}{*}{} & 100	&\textbf{0.735297} &	0.697297&	0.724628&	0.638984 \\ \cline{2-7}
				
				\multirow{9}{*}{ }&\multirow{3}{*}{8} & 20 &	\textbf{0.799082}	&	0.680402&	0.736225&	0.594415  \\ \cline{3-7}
				\multirow{9}{*}{}&\multirow{3}{*}{} & 50 &0.708704	&0.689910&	\textbf{0.711009}	&0.619342 \\ \cline{3-7}
				\multirow{9}{*}{}&\multirow{3}{*}{} & 100	&\textbf{0.781145} &	0.698146&	0.750615&	0.595450 \\ \hline
				
				\multirow{9}{*}{4}&\multirow{3}{*}{2} & 20 &	\textbf{0.758691}	&	0.689227&	0.733954&	0.624130  \\ \cline{3-7}
				\multirow{9}{*}{}&\multirow{3}{*}{} & 50	&\textbf{0.763147} &		0.690873&	0.748324&	0.637060 \\ \cline{3-7}
				\multirow{9}{*}{}&\multirow{3}{*}{} & 100	&\textbf{0.730720} &		0.713162 &	0.715489 &	0.640808 \\ \cline{2-7}
				
				\multirow{9}{*}{ }&\multirow{3}{*}{5} & 20 &	\textbf{0.801364}	 &		0.679019&	0.764121&	0.6955272  \\ \cline{3-7}
				\multirow{9}{*}{}&\multirow{3}{*}{} & 50	&\textbf{0.712999} &	0.692416	&0.684561&	0.6048977 \\ \cline{3-7}
				\multirow{9}{*}{}&\multirow{3}{*}{} & 100	&\textbf{0.775285} &	0.664577&	0.741605&	0.595998 \\ \cline{2-7}
				
				\multirow{9}{*}{ }&\multirow{3}{*}{8} & 20 &	\textbf{0.783354}	&		0.663554&	0.759846&	0.577628  \\ \cline{3-7}
				\multirow{9}{*}{}&\multirow{3}{*}{} & 50 &	\textbf{0.727478}	&	0.709145&	0.711456&	0.663443 \\ \cline{3-7}
				\multirow{9}{*}{}&\multirow{3}{*}{} & 100	&\textbf{0.758634} &	0.690469&	0.733429&	0.591219 \\ \hline

				\multirow{9}{*}{5}&\multirow{3}{*}{2} & 20 	&0.690899&	0.699662&	\textbf{0.701525}&		0.652180  \\ \cline{3-7}
				\multirow{9}{*}{}&\multirow{3}{*}{} & 50	&\textbf{0.713242} &	0.712468&	0.710150&	0.634902 \\ \cline{3-7}
				\multirow{9}{*}{}&\multirow{3}{*}{} & 100	 &0.698824&\textbf{0.716379}	&	0.704651&	0.648044 \\ \cline{2-7}
				
				\multirow{9}{*}{ }&\multirow{3}{*}{5} & 20 &	\textbf{0.749176}	 &	0.685902&	0.724796&	0.675969  \\ \cline{3-7}
				\multirow{9}{*}{}&\multirow{3}{*}{} & 50	&\textbf{0.729309} &	0.694974&	0.713021&	0.677386 \\ \cline{3-7}
				\multirow{9}{*}{}&\multirow{3}{*}{} & 100	&\textbf{0.763138} &	0.687546&	0.734786&	0.681083 \\ \cline{2-7}
				
				\multirow{9}{*}{ }&\multirow{3}{*}{8} & 20 &	\textbf{0.779347}	&	0.707880&	0.754157&	0.579095  \\ \cline{3-7}
				\multirow{9}{*}{}&\multirow{3}{*}{} & 50 &	\textbf{0.765622}	&	0.667469&	0.755479&	0.6015463 \\ \cline{3-7}
				\multirow{9}{*}{}&\multirow{3}{*}{} & 100	&\textbf{0.733192} &	0.710643&	0.715427&	0.637003 \\ \hline

				\multirow{9}{*}{6}&\multirow{3}{*}{2} & 20 &0.719806	&0.722424	&	\textbf{0.734201}&	0.647811  \\ \cline{3-7}
				\multirow{9}{*}{}&\multirow{3}{*}{} & 50	&\textbf{0.702079} &	0.693147&	0.710427&	0.6165592 \\ \cline{3-7}
				\multirow{9}{*}{}&\multirow{3}{*}{} & 100	 &\textbf{0.763388}	&	0.686958&	0.729897&	0.621244 \\ \cline{2-7}
				
				\multirow{9}{*}{ }&\multirow{3}{*}{5} & 20 &	\textbf{0.763755}	 &	0.702172&	0.734510&	0.693951  \\ \cline{3-7}
				\multirow{9}{*}{}&\multirow{3}{*}{} & 50	&\textbf{0.733391} &	0.686126&	0.719880&	0.630869 \\ \cline{3-7}
				\multirow{9}{*}{}&\multirow{3}{*}{} & 100	&\textbf{0.748348} &	0.712714&	0.722406&	0.5738283 \\ \cline{2-7}
				
				\multirow{9}{*}{ }&\multirow{3}{*}{8} & 20 &	\textbf{0.737870}	&	0.679402	&0.715045	&0.637772  \\ \cline{3-7}
				\multirow{9}{*}{}&\multirow{3}{*}{} & 50 &	\textbf{0.749314}	&	0.689410	&0.730157	&0.648392 \\ \cline{3-7}
				\multirow{9}{*}{}&\multirow{3}{*}{} & 100	&\textbf{0.774391} &	0.690839&	0.745160&	0.604552 \\ \hline
				
				\multicolumn{3}{|c|}{$Mean$} &\textbf{0.750615}	&0.696743&	0.728522&	0.633658  \\ \hline
							
			\end{tabular}
\end{table*}

\textcolor[rgb]{0,0,0}{It is obvious from Table \ref{Comparison of average $SP$  values} that KDMA achieves the best $SP$ values in forty-one instances, NSGA-II in one instance, MOEA/D in three instances, and WASFGA does not achieve any optimal $SP$ values  in all forty-five $\left ( f,m,n \right )$ combinations.} The mean of KDMA’s $SP$ values across the forty-five datasets is 0.750615 (better than those realized using the other three algorithms). Therefore, in terms of solving the EEDHFSSP, KDMA demonstrates superior solution diversity.

\begin{table*}[!ht]
	\center
	\caption{Comparison of average $GD$  values}
	\label{Comparison of average  $GD$ values}
		\setlength{\tabcolsep}{0.8cm}
			\begin{tabular}{|c|c|c|c|c|c|c|} \hline
			$F$ & $m$ & $n$ & KDMA & NSGA-II & MOEA/D & WASFGA  \\ \hline
				
				\multirow{9}{*}{2}&\multirow{3}{*}{2} & 20 &	\textbf{0.130907}	&	0.209836&	0.134251&	0.263656 \\ \cline{3-7}
				\multirow{9}{*}{}&\multirow{3}{*}{} & 50	&\textbf{0.117021} &	0.188481&	0.145242&	0.245695 \\ \cline{3-7}
				\multirow{9}{*}{}&\multirow{3}{*}{} & 100	&\textbf{0.131465} &	0.195991&	0.146456	&0.283036 \\ \cline{2-7}
				
				\multirow{9}{*}{ }&\multirow{3}{*}{5} & 20 &	\textbf{0.142824}	&	0.180450&	0.153150	&0.289138  \\ \cline{3-7}
				\multirow{9}{*}{}&\multirow{3}{*}{} & 50	 &0.139271	&0.185475&	\textbf{0.124258}&	0.239942 \\ \cline{3-7}
				\multirow{9}{*}{}&\multirow{3}{*}{} & 100	&\textbf{0.125911} & 	0.185522 &	0.133522& 	0.269129  \\ \cline{2-7}
				
				\multirow{9}{*}{ }&\multirow{3}{*}{8} & 20 &0.153784 	&0.192498 &\textbf{0.150404}	&0.271157   \\ \cline{3-7}
				\multirow{9}{*}{}&\multirow{3}{*}{} & 50	&\textbf{0.157753} & 	0.173053 &	0.159463 &	0.306227  \\ \cline{3-7}
				\multirow{9}{*}{}&\multirow{3}{*}{} & 100	&\textbf{0.166116 } &	0.183362 &	0.175642 &	0.245474  \\ \hline
				
				\multirow{9}{*}{3}&\multirow{3}{*}{2} & 20 &	\textbf{0.084211}	& 	0.188526 &	0.135264 &	0.208595   \\ \cline{3-7}
				\multirow{9}{*}{}&\multirow{3}{*}{} & 50	&\textbf{0.125353 } &		0.159955 &	0.126843 &	0.249466  \\ \cline{3-7}
				\multirow{9}{*}{}&\multirow{3}{*}{} & 100	&\textbf{0.087903} & 	0.208731 &	0.142685 &	0.275144  \\ \cline{2-7}
				
				\multirow{9}{*}{ }&\multirow{3}{*}{5} & 20 &	\textbf{0.104572}	 &	0.168560 &	0.122024 &	0.231684   \\ \cline{3-7}
				\multirow{9}{*}{}&\multirow{3}{*}{} & 50	&\textbf{0.149928} & 0.157242& 	0.151203 &	0.248649  \\ \cline{3-7}
				\multirow{9}{*}{}&\multirow{3}{*}{} & 100	&\textbf{0.124900} & 0.198360& 	0.149855& 	0.229609 \\ \cline{2-7}
				
				\multirow{9}{*}{ }&\multirow{3}{*}{8} & 20 &0.172210 	&0.174276 	&\textbf{0.156342}& 	0.238206   \\ \cline{3-7}
				\multirow{9}{*}{}&\multirow{3}{*}{} & 50 &\textbf{0.171589} 	&0.197791 &	0.178925 	&0.285797  \\ \cline{3-7}
				\multirow{9}{*}{}&\multirow{3}{*}{} & 100	&\textbf{0.150863} & 	0.178183 & 	0.170124 & 	0.226481  \\ \hline
				
				\multirow{9}{*}{4}&\multirow{3}{*}{2} & 20 &	\textbf{0.105579}	& 	0.180463 &	0.135825 &	0.270665   \\ \cline{3-7}
				\multirow{9}{*}{}&\multirow{3}{*}{} & 50	&\textbf{0.122612} & 	0.172362 & 	0.140457 & 	0.238674  \\ \cline{3-7}
				\multirow{9}{*}{}&\multirow{3}{*}{} & 100	&\textbf{0.094088} & 	0.241235& 	0.124646 &	0.315502  \\ \cline{2-7}
				
				\multirow{9}{*}{ }&\multirow{3}{*}{5} & 20 &	\textbf{0.149340}	 & 	0.165323  &	0.155462  &	0.221736   \\ \cline{3-7}
				\multirow{9}{*}{}&\multirow{3}{*}{} & 50	&\textbf{0.142910} & 	0.152168& 	0.151342& 	0.207915  \\ \cline{3-7}
				\multirow{9}{*}{}&\multirow{3}{*}{} & 100	&0.173359 	&\textbf{0.149749 }	&0.157421 	&0.265105  \\ \cline{2-7}
				
				\multirow{9}{*}{ }&\multirow{3}{*}{8} & 20 &0.169032 	&0.169018 	&\textbf{0.163450}& 	0.229791   \\ \cline{3-7}
				\multirow{9}{*}{}&\multirow{3}{*}{} & 50 &	\textbf{0.177563}	& 	0.188397& 	0.179642 &	0.301976  \\ \cline{3-7}
				\multirow{9}{*}{}&\multirow{3}{*}{} & 100	&\textbf{0.147685} & 	0.176995 &	0.150125 &	0.215146  \\ \hline

				\multirow{9}{*}{5}&\multirow{3}{*}{2} & 20 	&\textbf{0.073698} &	0.184855 	&0.135406 &	0.314665   \\ \cline{3-7}
				\multirow{9}{*}{}&\multirow{3}{*}{} & 50	&\textbf{0.143801} & 	0.174189 &	0.147982 &	0.210477  \\ \cline{3-7}
				\multirow{9}{*}{}&\multirow{3}{*}{} & 100	 &\textbf{0.105752} &0.200810 	&0.124568 &	0.272070  \\ \cline{2-7}
				
				\multirow{9}{*}{ }&\multirow{3}{*}{5} & 20 &	\textbf{0.113409}	 & 	0.150964 &	0.126431 &	0.208673   \\ \cline{3-7}
				\multirow{9}{*}{}&\multirow{3}{*}{} & 50	&\textbf{0.159253} & 	0.169624&  	0.160158 & 	0.263108  \\ \cline{3-7}
				\multirow{9}{*}{}&\multirow{3}{*}{} & 100	&\textbf{0.154368} & 	0.159390 &	0.160013 &	0.244988  \\ \cline{2-7}
				
				\multirow{9}{*}{ }&\multirow{3}{*}{8} & 20 &	\textbf{0.152244}	& 	0.153972 &	0.155423 &	0.248254   \\ \cline{3-7}
				\multirow{9}{*}{}&\multirow{3}{*}{} & 50 &	\textbf{0.147033}	& 	0.201189 &	0.164215 &	0.256551  \\ \cline{3-7}
				\multirow{9}{*}{}&\multirow{3}{*}{} & 100	&0.168920 	&0.175371 &\textbf{0.166478}	 &	0.267107  \\ \hline

				\multirow{9}{*}{6}&\multirow{3}{*}{2} & 20 &\textbf{0.088489} 	&0.183778 	&0.135264 	&0.334893   \\ \cline{3-7}
				\multirow{9}{*}{}&\multirow{3}{*}{} & 50	&\textbf{0.118600} & 	0.150759  &	0.142106  &	0.202783  \\ \cline{3-7}
				\multirow{9}{*}{}&\multirow{3}{*}{} & 100	 &\textbf{0.175823}	& 	0.180070 &	0.179543 &	0.272389 \\ \cline{2-7}
				
				\multirow{9}{*}{ }&\multirow{3}{*}{5} & 20 &	0.162428	 & \textbf{0.146049}	& 	0.154217 &	0.190410   \\ \cline{3-7}
				\multirow{9}{*}{}&\multirow{3}{*}{} & 50	&\textbf{0.163159} & 	0.183072 &	0.175465 &	0.212960  \\ \cline{3-7}
				\multirow{9}{*}{}&\multirow{3}{*}{} & 100	&\textbf{0.134627} & 	0.182143 &	0.143247 &	0.224290  \\ \cline{2-7}
				
				\multirow{9}{*}{ }&\multirow{3}{*}{8} & 20 &	\textbf{0.160114}	& 	0.207423 &	0.169746 &	0.279497   \\ \cline{3-7}
				\multirow{9}{*}{}&\multirow{3}{*}{} & 50 &	\textbf{0.154402}	& 	0.176949 &	0.166485 &	0.320717  \\ \cline{3-7}
				\multirow{9}{*}{}&\multirow{3}{*}{} & 100	&\textbf{0.165445} & 	0.186427& 	0.170124 &	0.203172 \\ \hline
				
				\multicolumn{3}{|c|}{$Mean$} &\textbf{0.139118}	&	0.179756&	0.150908&	0.253346 \\ \hline
							
			\end{tabular}
\end{table*}

According to Table \ref{Comparison of average  $GD$ values}, KDMA has the best average values in thirty-eight instances, NSGA-II in two instances, MOEA/D in five instances, and WASFGA in none. Additionally, the overall average value (0.139118) of KDMA is significantly lower than those of NSGA-II, MOEA/D, and WASFGA, indicating the superior convergence of KDMA.

\begin{table*}[!ht]
	\center
	\caption{Comparison of average $IGD$ values}
	\label{Comparison of average $IGD$ values}
		\setlength{\tabcolsep}{0.8cm}
			\begin{tabular}{|c|c|c|c|c|c|c|} \hline
			$F$ & $m$ & $n$ & KDMA & NSGA-II & MOEA/D & WASFGA  \\ \hline
				
				\multirow{9}{*}{2}&\multirow{3}{*}{2} & 20 &	\textbf{0.099725}	&	0.165646&	0.153325&	0.194498 \\ \cline{3-7}
				\multirow{9}{*}{}&\multirow{3}{*}{} & 50	&\textbf{0.126656} &	0.163129&	0.144145&	0.243189 \\ \cline{3-7}
				\multirow{9}{*}{}&\multirow{3}{*}{} & 100	&\textbf{0.121372} &	0.183355&	0.185146	&0.220639\\ \cline{2-7}
				
				\multirow{9}{*}{ }&\multirow{3}{*}{5} & 20 &0.145625 	&0.167178 	&\textbf{0.122575} &	0.140707   \\ \cline{3-7}
				\multirow{9}{*}{}&\multirow{3}{*}{} & 50	 &\textbf{0.137474} 	&0.190675 	&0.187606& 	0.211841  \\ \cline{3-7}
				\multirow{9}{*}{}&\multirow{3}{*}{} & 100	&\textbf{0.149501} & 	0.155092 &	0.158336 &	0.186242  \\ \cline{2-7}
				
				\multirow{9}{*}{ }&\multirow{3}{*}{8} & 20 &\textbf{0.129004} 	&0.158984 	&0.146007 	&0.205527    \\ \cline{3-7}
				\multirow{9}{*}{}&\multirow{3}{*}{} & 50	&\textbf{0.144635} & 	0.191440 &	0.194440 &	0.247086  \\ \cline{3-7}
				\multirow{9}{*}{}&\multirow{3}{*}{} & 100	&\textbf{0.125820} & 	0.143361 &	0.141841 &	0.181570  \\ \hline
				
				\multirow{9}{*}{3}&\multirow{3}{*}{2} & 20 &	\textbf{0.098842}	& 	0.155548 &	0.122575 &	0.225182 \\ \cline{3-7}
				\multirow{9}{*}{}&\multirow{3}{*}{} & 50	&\textbf{0.108451} & 	0.189477 &	0.129891 &	0.203021  \\ \cline{3-7}
				\multirow{9}{*}{}&\multirow{3}{*}{} & 100	&\textbf{0.071217} & 	0.208990 &	0.142015 &	0.187182   \\ \cline{2-7}
				
				\multirow{9}{*}{ }&\multirow{3}{*}{5} & 20 &	\textbf{0.120471}	 & 	0.149720  &	0.135422  &	0.163979 \\ \cline{3-7}
				\multirow{9}{*}{}&\multirow{3}{*}{} & 50	&\textbf{0.139077 } &	0.155620 &	0.151055 	&0.178884  \\ \cline{3-7}
				\multirow{9}{*}{}&\multirow{3}{*}{} & 100	&\textbf{0.104930} &  	0.176947 &	0.161120 &	0.201196 \\ \cline{2-7}
				
				\multirow{9}{*}{ }&\multirow{3}{*}{8} & 20 &0.173000 	&0.165429 &\textbf{0.154284}	 &	0.179023    \\ \cline{3-7}
				\multirow{9}{*}{}&\multirow{3}{*}{} & 50 &\textbf{0.136784} 	& 	0.167173 &	0.147855 &	0.166502  \\ \cline{3-7}
				\multirow{9}{*}{}&\multirow{3}{*}{} & 100	&\textbf{0.153265} & 	0.140916 & 	0.153421 & \textbf{0.139563}	   \\ \hline
				
				\multirow{9}{*}{4}&\multirow{3}{*}{2} & 20 &	\textbf{0.087671}	& 	0.173094& 	0.124251 	&0.186692   \\ \cline{3-7}
				\multirow{9}{*}{}&\multirow{3}{*}{} & 50	&\textbf{0.105270} & 	0.163451 &	0.135250 &	0.198447  \\ \cline{3-7}
				\multirow{9}{*}{}&\multirow{3}{*}{} & 100	&\textbf{0.087575} &  	0.154939 	&0.142121 &	0.170695  \\ \cline{2-7}
				
				\multirow{9}{*}{ }&\multirow{3}{*}{5} & 20 &0.163019 &\textbf{0.136626}	 &	0.142511 	&0.169446    \\ \cline{3-7}
				\multirow{9}{*}{}&\multirow{3}{*}{} & 50	&0.150840 &	0.153240 &\textbf{0.133528}	 	&0.144785   \\ \cline{3-7}
				\multirow{9}{*}{}&\multirow{3}{*}{} & 100	&\textbf{0.153136} & 0.169468 &	0.166454 	&0.168217   \\ \cline{2-7}
				
				\multirow{9}{*}{ }&\multirow{3}{*}{8} & 20 &0.164111 &\textbf{0.144535}	 &	0.152435 &	0.157462    \\ \cline{3-7}
				\multirow{9}{*}{}&\multirow{3}{*}{} & 50 &	0.147515 	&\textbf{0.144261} &	0.147850 &	0.151218   \\ \cline{3-7}
				\multirow{9}{*}{}&\multirow{3}{*}{} & 100	&\textbf{0.144501} & 	0.156920  &	0.155648 & 	0.163399   \\ \hline

				\multirow{9}{*}{5}&\multirow{3}{*}{2} & 20 	&\textbf{0.061719} & 	0.166276 &	0.138521 &	0.252240   \\ \cline{3-7}
				\multirow{9}{*}{}&\multirow{3}{*}{} & 50	&\textbf{0.135102} & 	0.170598 &	0.158498 &	0.184704  \\ \cline{3-7}
				\multirow{9}{*}{}&\multirow{3}{*}{} & 100	 &\textbf{0.101984} & 	0.178423 &	0.169967 &	0.226277  \\ \cline{2-7}
				
				\multirow{9}{*}{ }&\multirow{3}{*}{5} & 20 &	\textbf{0.100546}	 &	0.138746 &	0.129785 &	0.151350   \\ \cline{3-7}
				\multirow{9}{*}{}&\multirow{3}{*}{} & 50	&0.139545 	&0.138314 &\textbf{0.136354}	 &	0.162086   \\ \cline{3-7}
				\multirow{9}{*}{}&\multirow{3}{*}{} & 100	&\textbf{0.142699} & 	0.172142 &	0.155846 &	0.176521  \\ \cline{2-7}
				
				\multirow{9}{*}{ }&\multirow{3}{*}{8} & 20 &	\textbf{0.136523}	& 	0.152567 &	0.164010 	&0.193828   \\ \cline{3-7}
				\multirow{9}{*}{}&\multirow{3}{*}{} & 50 &	\textbf{0.104025}	& 	0.145856& 	0.151785 &	0.196455   \\ \cline{3-7}
				\multirow{9}{*}{}&\multirow{3}{*}{} & 100	&\textbf{0.134346} &	0.156988 &	0.142555 	&0.186240   \\ \hline

				\multirow{9}{*}{6}&\multirow{3}{*}{2} & 20 &\textbf{0.084394} 	& 	0.162152 & 	0.158421 & 	0.240563   \\ \cline{3-7}
				\multirow{9}{*}{}&\multirow{3}{*}{} & 50	&\textbf{0.102837} & 	0.140105 &	0.141425 	&0.159635  \\ \cline{3-7}
				\multirow{9}{*}{}&\multirow{3}{*}{} & 100	 &\textbf{0.151013}	& 	0.173533 &	0.163565 &	0.221553  \\ \cline{2-7}
				
				\multirow{9}{*}{ }&\multirow{3}{*}{5} & 20 &\textbf{0.118012} &	0.134985 &	0.144882 &	0.147451   \\ \cline{3-7}
				\multirow{9}{*}{}&\multirow{3}{*}{} & 50	&\textbf{0.1131699} & 	0.156111  & 	0.150452  & 	0.157785   \\ \cline{3-7}
				\multirow{9}{*}{}&\multirow{3}{*}{} & 100	&\textbf{0.1223997} & 	0.158768 &	0.133649 &	0.155864   \\ \cline{2-7}
				
				\multirow{9}{*}{ }&\multirow{3}{*}{8} & 20 &	\textbf{0.125170}	& 	0.163610 &	0.145206 	&0.243423   \\ \cline{3-7}
				\multirow{9}{*}{}&\multirow{3}{*}{} & 50 &	\textbf{0.128028}	& 	0.170973 &	0.153974 	&0.210269  \\ \cline{3-7}
				\multirow{9}{*}{}&\multirow{3}{*}{} & 100	&\textbf{0.119191} & 	0.137920 &	0.131022 	&0.156111  \\ \hline
				
				\multicolumn{3}{|c|}{$Mean$} &\textbf{0.1246718}	&	0.160962	&    0.144891&	0.186856 \\ \hline
							
			\end{tabular}
\end{table*}

\textcolor[rgb]{0,0,0}{The experimental results of \textit{IGD} are shown in Table \ref{Comparison of average $IGD$ values}. From Table \ref{Comparison of average $IGD$ values}, it can be seen that KDMA has the best performance. Regarding the overall average $IGD$, the average value (0.124671) of KDMA is significantly lower than those of the other three algorithms. Hence, KDMA demonstrates superior overall performance in solving the EEDHFSSP and holds a clear advantage, especially when applied to large-scale datasets.}

\textcolor[rgb]{0,0,0}{In summary, KDMA surpasses the comparison algorithms in multiple evaluation aspects. Moreover, certain strategies employed in KDMA can also be applied to address energy-efficient concerns in manufacturing system scheduling problems. The remarkable performance of KDMA in solving EEDHFSSP demonstrates the algorithm's significant practical value. Certainly, our efforts to enhance the algorithm are merely preliminary attempts, and there are still numerous directions that warrant further exploration. For example, machine learning are effective methods to complex combinatorial optimization problems such as EEDHFSSP \cite{song2023ensemble}.}

\section{Conclusion}
\label{Conclusion}
\textcolor[rgb]{0,0,0}{With the growing global attention to climate change, research and applications related to carbon emissions, emission peaks, and carbon neutrality are gradually unfolding in the field of production scheduling. In this study, we designed scheduling model and solution algorithms for EEDHFSSP. Initially, a mathematical model of EEDHFSSP is established with the optimization objectives of minimizing both the makespan and total carbon emissions, and the methods for calculating completion time and total carbon emissions are provided. Then, KDMA is proposed to solve EEDHFSSP. To improve the quality of initial solutions, a collaborative initialization strategy is introduced. New update operations are proposed to strengthen global search capabilities, and a local search strategy based on key factories is developed to improve algorithm performance. Moreover, a carbon reduction strategy is proposed to decrease carbon emissions. Finally, experiments are conducted on problem-specific instances. Experimental results show that KDMA is more advantageous for solving EEDHFSP compared to other algorithms.}

\textcolor[rgb]{0,0,0}{In future research, we will consider other DHFSSPs that closely resemble real-world scenarios, such as equipment failures or temporary postponement or cancellation of jobs. These intricate situations will impose higher demands on the algorithms. It is worth exploring machine learning and reinforcement learning methods to facilitate algorithm evolution. Effective problem feature mining can enhance the ability of optimization algorithms to discover high-quality solutions.}

\section*{Declaration of Competing Interest}
The authors declare that they have no known competing financial interests or personal relationships that could have appeared to influence the work reported in this paper.


\bibliographystyle{unsrt}

\bibliography{mybib}

\end{document}